\documentclass[a4paper,11pt]{amsart}

\usepackage{amssymb,epsf,verbatim,hyperref,showkeys}

\theoremstyle{plain}
\newtheorem{lem}{Lemma}[section]
\newtheorem{prop}[lem]{Proposition}
\newtheorem{cor}[lem]{Corollary}
\newtheorem{thm}[lem]{Theorem}
\theoremstyle{definition}
\newtheorem{Def}[lem]{Definition}
\theoremstyle{remark}
\newtheorem*{note}{Note}
\newtheorem*{eg*}{Example}

\newtheorem*{rem}{Remark}

\newcommand{\func}{\operatorname}
\newcommand{\up}[1]{\text{$\uparrow^{#1}$}}
\newcommand{\down}[1]{\text{$\downarrow^{#1}$}}

\begin{document}

\title{Beyond Rouquier partitions}
\author{Kai Meng Tan}
\address{Department of Mathematics, National University of
Singapore, 2, Science Drive 2, Singapore 117543.}
\email{tankm@nus.edu.sg}
\date{Sep 2006}
\thanks{Supported by Academic Research Fund R-146-000-043-112 of
National University of Singapore.}
\subjclass[2000]{Primary: 20C08; Secondary: 17B37}
%\keywords{Fock space, canonical bases, decomposition numbers, Hecke
%algebras, quantized Schur algebras}

\begin{abstract}
We obtain closed formulas, in terms of Littlewood-Richardson coefficients, for the canonical basis elements of the Fock space representation of $U_v(\widehat{\mathfrak{sl}}_e)$ which are labelled by partitions having `locally small' $e$-quotients and arbitrary $e$-cores.  We further show that upon evaluation at $v=1$, this gives the corresponding decomposition numbers of the $q$-Schur algebra in characteristic $l$ (where $q$ is a primitive $e$-th root of unity if $l \ne e$ and $q=1$ otherwise) whenever $l$ is greater than the size of each constituent of the $e$-quotient.
\end{abstract}
\maketitle

\section{Introduction}

In the course of investigating the truth of Brou\'{e}'s Abelian defect group conjecture for symmetric groups, Rouquier \cite{Rouquier-Hab} singled out a special class of blocks of symmetric groups which he believed to have good properties.  These blocks, and their corresponding blocks of the Iwahori-Hecke algebras (of type A) and the $q$-Schur algebras, are now known as Rouquier blocks, and are well understood, by the works of several authors (see, for example, \cite{CK}, \cite{mr}, \cite{Miyachi-thesis}, \cite{Turner}, \cite{JML}).  In particular, there exist closed formulas, in terms of Littlewood-Richardson coefficients, for the decomposition numbers of these blocks when they are of `Abelian defect'.

These formulas coincide with those for the $v$-decomposition numbers arising from the canonical basis of the Fock space representation of $U_v(\widehat{\mathfrak{sl}}_e)$  obtained by Leclerc and Miyachi \cite{LM} upon evaluation at $v=1$.  In fact, such formulas for the $v$-decomposition numbers hold not only for the Rouquier partitions --- the partitions indexing the simple modules of Rouquier blocks.  As shown by Chuang and the author \cite{can}, whenever $\kappa$ is an $e$-core partition of {\em Rouquier type} --- one having an abacus display in which the number of beads on each runner is non-decreasing as we go from left to right --- then the formulas hold for the canonical basis element labelled by any partition lying in the class $\mathcal{P}_{\kappa}$ of $e$-regular partitions having $e$-core $\kappa$ and `locally small' $e$-quotients.

However, this larger class of partitions is unsatisfactory for the following reasons:
\begin{itemize}
\item It does not contain any $e$-singular partition.
\item It is not closed under the Mullineux involution $\mu \mapsto \mu^*$; i.e.\ $\mu \in \mathcal{P}_{\kappa}$ does not imply that $\mu^* \in \mathcal{P}_{\kappa'}$ (even though, as shown in Proposition \ref{P:Mull}, that if the formula holds for the canonical basis element labelled by $\mu$, then it also holds for that labelled by $\mu^*$).
\item $\kappa$ is not arbitrary.
\end{itemize}
We address all the above shortcomings in the first part of this paper.  Our class $\mathcal{P}_{\kappa}^*$, consisting of partitions having $e$-core $\kappa$ and `locally small' $e$-quotients (though our definition of `locally small' is different from that for $\mathcal{P}_{\kappa}$), is defined for all $e$-core partitions $\kappa$, is closed under the Mullineux involution, and consists of $e$-singular partitions as well.  Furthermore, if $\kappa$ is of Rouquier type, then $\mathcal{P}_{\kappa}$ is a proper subclass of $\mathcal{P}_{\kappa}^*$.

In the second part of this paper, we show that, upon evaluation at $v=1$, our formula for the canonical basis element labelled by a partition $\mu \in \mathcal{P}_{\kappa}^*$ gives the corresponding decomposition numbers of the $q$-Schur algebras in characteristic $l$ (where $q$ is a primitive $e$-th root of unity if $l \ne e$, and $q=1$ otherwise) whenever $l$ is larger than the sizes of each constituent of the $e$-quotient of $\mu$.  This thus includes cases where the blocks, in which the simple modules labelled by such partitions lie, are not of `Abelian defect'.

The paper is organized as follows: we begin in section 2 with a short account of the background theory which we require.  In section 3, we introduce the notations used in this paper, and state the main theorem.  Sections 4 and 5 are devoted to the proof of the main theorem.

\section{Preliminaries}

In this section, we give a brief account of the background theory we require, and introduce some notations and conventions which will be used in this paper.  From now on, we fix an integer $e \geq 2$.

\subsection{Partitions}
Let $\mathcal{P}_{n}$ be the set of partitions of a natural number $n$, and denote the lexicographical ordering and dominance ordering on $\mathcal{P}_n$ by $\geq$ and $\trianglerighteq$ respectively.  Write $\mathcal{P}=\bigcup_{n} \mathcal{P}_n$ for the set of all partitions.

Let $\lambda =(\lambda _{1},\lambda _{2},\dotsc ,\lambda _{s})$ be a
partition, with $\lambda _{1}\geq \lambda _{2}\geq \dotsb \geq \lambda
_{s}>0 $. We write $|\lambda |=\sum_{i=1}^{s}\lambda _{i}$ and $l(\lambda
)=s $.  The conjugate partition of $\lambda$ will be denoted as $\lambda'$.

The James $e$-abacus (see, for example,
\cite[Section 2.7]{JK}) has $e$ runners, numbered $0$ to $e-1$ from left to right, and its positions are numbered from left to right and down the rows, starting from $0$.  Any partition may be displayed on such an abacus; this display is only unique when the number of
beads used is fixed. In this paper, if $\lambda $ is a partition with $e$-core $\kappa $, we shall always display $\lambda$ on an $e$-abacus with $l(\kappa) + Ne$ beads for sufficiently large $N$, and denote its $e$-quotient by $(\lambda ^{0},\lambda ^{1},\dotsc
,\lambda ^{n-1})$, where $\lambda ^{i}$ is the partition read off from runner $i$ of the abacus display of $\lambda $.  Define the $e$-weight of a bead in the abacus display of $\lambda$ to be the number of vacant positions above it in the same runner.  The $e$-weight of a runner is the sum of $e$-weights of all the beads in the runner, and the $e$-weight of $\lambda$ is the sum of $e$-weights of all the beads in the abacus display of $\lambda$.

If $s+t=r$ and $\alpha \in \mathcal{P}_s$ and $\beta \in \mathcal{P}_t$, let $c_{\alpha \beta }^{\lambda }$ denote
the multiplicity of the ordinary irreducible character $\chi ^{\lambda }$ of the symmetric group $\mathfrak{S}_r$ in the induced character $\func{Ind}_{\mathfrak{S}%
_{s}\times \mathfrak{S}_{t}}^{\mathfrak{S}_{r}}(\chi ^{\alpha }\otimes \chi
^{\beta })$. We refer the reader to \cite[I.9]{Mac} for a combinatorial
description of $c_{\alpha \beta }^{\lambda }$, which is known as a
Littlewood-Richardson coefficient. By convention, we define $c_{\alpha \beta
}^{\lambda }$ to be $0$ when $|\lambda |\neq |\alpha |+| \beta| $.

\subsection{The Fock space representation}

The algebra $U_{v}(\widehat{\mathfrak{sl}}_{e})$ is the
associative algebra over $\mathbb{C}(v)$ with generators $e_{r}$,
$f_{r}$, $k_{r}$, $k_{r}^{-1}$ $(0\leq r\leq e-1)$, $d$, $d^{-1}$
subject to some relations (see, for example, \cite[\S4]{Leclerc2}).
An important $U_{v}(\widehat{\mathfrak{sl}}_{e})$-module is the Fock space representation
$\mathcal{F}$, which as a $\mathbb{C}(v)$-vector space has $\mathcal{P}$ as a basis.  For our purposes, an explicit description of the actions of $e_r$ and $f_r$ will suffice.

Display a partition $\lambda$ on the $e$-abacus with $t$ beads, where $t \geq l(\lambda)$ and $e \nmid (r+t)$.  Let $i$ be the residue class of $(r+t)$ modulo $e$.  Suppose there is a bead on runner $i-1$ whose succeeding position on runner $i$ is vacant; let $\mu$ be the partition obtained when this bead is moved to its succeeding position.
Let $N_>(\lambda,\mu)$ (resp.\ $N_<(\lambda,\mu)$) be the number of beads on runner $i-1$ below (resp.\ above) the bead moved to obtained $\mu$ minus the number of beads on runner $i$ below (resp.\ above) the vacant position that becomes occupied in obtaining $\mu$.
We have
\begin{align*}
f_r( s(\lambda)) &= \sum_{\mu} v^{N_>(\lambda,\mu)} s(\mu); \\
e_r( s(\mu)) &= \sum_{\lambda} v^{-N_<(\lambda,\mu)} s(\lambda),
\end{align*}
where $\mu$ in the first sum runs over all partitions that can be obtained from $\lambda$ by moving a bead on runner $i-1$ to its vacant succeeding position on runner $i$, while $\lambda$ in the second sum runs over all partitions that can be obtained from $\mu$ by moving a bead on runner $i$ to its vacant preceding position on runner $i-1$.

In \cite{LT1}, Leclerc and Thibon introduced an involution $x \mapsto
\overline{x}$ on $\mathcal{F}$, having the following properties (among
others):
$$\overline{a(v)x } = a(v^{-1})\overline{x}, \quad
\overline{e_r(x)} = e_r(\overline{x}), \quad
\overline{f_r(x)} = f_r(\overline{x}) \quad (a \in \mathbb{C}(v),\ x
\in \mathcal{F}).$$
For $k \in \mathbb{Z}^+$, we write $[k]! = \prod_{i=1}^k \frac{v^i - v^{-i}}{v-v^{-1}}$, and let $e_r^{(k)} = e_r^k/[k]!$ and $f_r^{(k)} = f_r^k/[k]!$.  Note that $\overline{e_r^{(k)}(x)} = e_r^{(k)}(\overline{x})$ and
$\overline{f_r^{(k)}(x)} = f_r^{(k)}(\overline{x})$.
There is a distinguished basis $\{G(\sigma ) \mid \sigma \in \mathcal{P}
\}$ of $\mathcal{F}$, called the canonical basis,
having the following
characterization (\cite[Theorem 4.1]{LT1}):
\begin{enumerate}
\item  $G(\sigma )- \sigma \in vL$, where $L$ is the
       $\mathbb{Z}[v]$-lattice in $\mathcal{F}$ generated by $\mathcal{P}$.
\item  $\overline{G(\sigma )}=G(\sigma )$.
\end{enumerate}
The $v$-decomposition number $d_{\lambda\sigma}(v)$ is defined the coefficient of $\lambda$
in $G(\sigma)$; these numbers enjoy the following property:

\begin{thm}[{\cite[Theorem 9, Proposition 11, Corollary 14]{Leclerc2}}] \label{LT1}
We have
\begin{align*}
d_{\sigma\sigma}(v) &= 1, \\
d_{\lambda\sigma}(v) &\in v \mathbb{N}_0[v] \ \text{ for all } \lambda \ne \sigma.
\end{align*}
Furthermore, $d_{\lambda\sigma}(v) \ne 0$ only if $\sigma
 \trianglerighteq \lambda$ and $\lambda$ and $\sigma$ have the same
 $e$-core.
\end{thm}

There is an involution $\mu \mapsto \mu^*$, known as the Mullineux involution, on the set of $e$-regular partitions (see, for example, \cite[page 120]{MathasBook}).  The $v$-decomposition numbers $d_{\lambda\mu}(v)$ and $d_{\lambda'\mu^*}(v)$ are related in the following way:

\begin{thm}[{\cite[Theorem 7.2]{LLT1}}] \label{T:Mull}
We have $d_{\lambda'\mu^*}(v) = v^w d_{\lambda\mu}(v^{-1})$ where $w$ is the $e$-weight of $\mu$.
\end{thm}

The partition ${\mu^*}'$ can be characterized in the following way using the $v$-decomposition numbers.

\begin{thm}[{\cite[Corollary 7.7]{LLT1}}] \label{T:Mull2}
We have $d_{{\mu^*}'\mu}(v) = v^w$, where $w$ is the $e$-weight of $\mu$, and $\deg d_{\lambda\mu}(v) < w$ for all $\lambda \ne {\mu^*}'$.
\end{thm}

\subsection{$q$-Schur algebras}

Let $\mathbb{F}$ be a field of characteristic $l$, and let $q \in \mathbb{F}^*$ be such that $e$ is the least integer such that $1+q+\dotsc+q^{e-1}= 0$.  The $q$-Schur algebra $\mathcal{S}_{\mathbb{F},q}(n) = \mathcal{S}_{\mathbb{F},q}(n,n)$ over $\mathbb{F}$ has a distinguished class $\{ \Delta^{\mu} \mid \mu \in \mathcal{P}_n \}$ of right modules called Weyl modules.  Each $\Delta^{\mu}$ has a simple head $L^{\mu}$, and the set $\{ L^{\mu} \mid \mu \in \mathcal{P}_n \}$ is a complete set of non-isomorphic simple modules of $\mathcal{S}_{\mathbb{F},q}(n)$.  The projective cover $P^\mu$ of $L^\mu$ (or of $\Delta^\mu$) has a filtration in which each factor is isomorphic to a Weyl module; the multiplicity of $\Delta^{\lambda}$ in such a filtration is well-defined, and is equal to the multiplicity of $L^{\mu}$ as a composition factor of $\Delta^\lambda$.  We denote this multiplicity as $d_{\lambda\mu}^l$, which is a decomposition number of $\mathcal{S}_{\mathbb{F},q}(n)$.

The decomposition numbers in characteristic $l$ and those in characteristic $0$ are related by an adjustment matrix $A_l$: let $D_l = (d_{\lambda\mu}^l)_{\lambda,\mu \in \mathcal{P}_n}$, then $D_l = D_0 A_l$.  Furthermore, the matrix $A_l$ is lower unitriangular with nonnegative entries when the partitions indexing its rows and columns are ordered by a total order extending the dominance order on $\mathcal{P}_n$ (such as the lexicographic order).  As a consequence, we have

\begin{lem}
$d_{\lambda\mu}^l \geq d_{\lambda\mu}^0$.
\end{lem}

The link between the $v$-decomposition numbers of the Fock space and the decomposition numbers of $q$-Schur algebras is established by Varagnolo and Vasserot:

\begin{thm}[\cite{VV}]
$d_{\lambda\mu}(1) = d_{\lambda\mu}^0$.
\end{thm}

Thus the canonical basis vector $G(\mu)$ of $\mathcal{F}$ corresponds to the projective cover $P^{\mu}$ of $q$-Schur algebras, while the standard basis element $\lambda$ of $\mathcal{F}$ corresponds to the Weyl module $\Delta^{\lambda}$.  Under this correspondence, the action of $e_r, f_r \in U_v(\widehat{\mathfrak{sl}}_e)$ on $\mathcal{F}$ corresponds to $r$-restriction and $r$-induction (see, for example, \cite[6.4]{MathasBook}) of modules of $q$-Schur algebras.

\subsection{Jantzen order}  Let $\lambda$ be a partition, and consider its abacus display, with $k$ beads say.  Suppose in moving a bead, say at position $a$, up its runner to some vacant position, say $a-ie$, we obtain (the abacus display of) a partition $\mu$.  Write $l_{\lambda \mu}$ for the number of occupied positions between $a$ and $a-ie$, and let $h_{\lambda\mu} = i$.  Also, write $\lambda \xrightarrow{\mu} \tau$ if the abacus display of $\mu$ with $k$ beads is also obtained from that of $\tau$ by moving a bead at position $b$ to $b-ie$, and $a < b$.  Thus if $\lambda \xrightarrow{\mu} \tau$, then the abacus display of $\tau$ with $k$ beads may be obtained from $\lambda$ in two steps: first move the bead at position $a$ to position $a-ie$ (which yields the abacus display of $\mu$), and then move the bead at position $b-ie$ to position $b$.

The Jantzen sum formula (see, for example, \cite[5.32]{MathasBook}) provides an upper bound for the decomposition numbers, and may be stated as follows:

\begin{thm} \label{T:Jantzen}
Let $\lambda$ and $\mu$ be distinct partitions, and let
$$J_{\lambda\mu} = \sum (-1)^{l_{\lambda\sigma} + l_{\tau\sigma} + 1} (1+\nu_l(h_{\lambda \sigma})) d_{\tau\mu}^l,$$
where the sum runs through all $\tau$ and $\sigma$ such that $\lambda \xrightarrow{\sigma} \tau$, and where $\nu_l$ denotes the standard $l$-valuation if $l > 0$ and $\nu_0(x) = 0$ for all $x$.
Then $d_{\lambda \mu}^l \leq J_{\lambda\mu}$, and $d_{\lambda \mu}^l = 0$ if and only if $J_{\lambda\mu} = 0$.
\end{thm}

We write $\lambda \rightarrow \tau$ if there exists some $\mu$ such that $\lambda \xrightarrow{\mu} \tau$.  We further write $\lambda <_J \sigma$ (or $\sigma >_J \lambda$) if there exist partitions $\tau_0,\dotsc, \tau_r$ such that $\tau_0 = \lambda$, $\tau_r = \sigma$ and $\tau_{i-1} \to \tau_i$ for all $i = 1,\dotsc, r$.
It is clear that $\geq_J$ (which means $>_J$ or $=$) defines a partial order on the set $\mathcal{P}$ of all partitions, and that if $\lambda \geq_J \mu$, then $\lambda$ and $\mu$ have the same $e$-core and $e$-weight.  Furthermore, the dominance order extends $\geq_J$.

The Jantzen sum formula motivates the study of $\leq_J$, as it provides this easy consequence.

\begin{lem} \label{L:J}
Suppose $d_{\lambda\mu}^l \ne 0$.  Then $\lambda \leq_J \mu$.
\end{lem}

Thus, if $\lambda \nleq_J \mu$, then $d_{\lambda\mu}^l = 0$.  Unfortunately, it is not easy to determine by inspection if $\lambda \nleq_J \mu$.  To this end, we introduce another partial order $\geq_p$.

Let $\lambda$ be a partition, and suppose that when displayed on an $e$-abacus with $N$ beads, the beads having positive $e$-weights are at positions $a_1,a_2,\dotsc, a_r$, with say the bead at position $a_i$ having $e$-weight $w_i$.  The {\em induced $e$-sequence} of $\lambda$, denoted $s(\lambda) = s(\lambda)_N$, is defined as
$$ \bigsqcup_{i=1}^r (a_i,a_i-e,\dotsc, a_i- (w_i-1)e),$$
where $(b_1,b_2,\dotsc, b_s) \sqcup (c_1,c_2,\dotsc,c_t)$ denotes the (unique) weakly decreasing sequence obtained by rearranging the terms in the sequence $(b_1,b_2,\dotsc, b_s,c_1,c_2,\dotsc,c_t)$.

\begin{eg*}
Let $\lambda = (6,6,5,4)$.  Its 3-abacus display with $6$ beads is
$$
\begin{smallmatrix}
\bullet & \bullet & - \\
- & - & - \\
\bullet & - & \bullet \\
- & \bullet & \bullet
\end{smallmatrix}.
$$
Thus $s(\lambda)_6 = (11,8) \sqcup (10,7) \sqcup (8,5) \sqcup (6) = (11,10,8,8,7,6,5)$.
\end{eg*}

We record some basic properties of induced $e$-sequences.

\begin{lem} \label{L:ez2}
Let $\lambda$ be a partition.  We have
\begin{enumerate}
\item $s(\lambda) \in \mathbb{N}_0^w$ where $w$ is the $e$-weight of $\lambda$;
\item if $M_j$ and $m_j$ denote respectively the largest occupied position and least vacant position on runner $j$ of the abacus display of $\lambda$, and $s(\lambda) = (l_1,l_2,\dotsc,l_w)$, then $M_j= \max_i\{ l_i \mid l_i \equiv j \pmod e \}$ and $m_j = \min_i \{ l_i \mid l_i \equiv j \pmod e \} - e$ whenever runner $j$ has positive $e$-weight.
\item if $\mu$ is the partition obtained by moving a bead in the abacus display of $\lambda$ (with $N$ beads) from position $x$ to position $x - ie$, then $s(\lambda)_N = s(\mu)_N \sqcup (x,x-e,\dotsc,x-(i-1)e)$.
\end{enumerate}
\end{lem}

The partial order $\geq_p$ on the set of partitions is defined as:  $\lambda \geq_p \mu$ if and only if $\lambda$ and $\mu$ have the same $e$-core and the same $e$-weight, say $w$, and $s(\lambda)_N\geq s(\mu)_N$ (for sufficiently large $N$) in the standard product order on $\mathbb{N}_0^w$.

\begin{lem} \label{L:Jimplyp}
If $\lambda \leq_J \mu$, then $\lambda\leq_p \mu$.
\end{lem}

\begin{proof}
It suffices to prove that $\lambda <_p \mu$ when $\lambda \xrightarrow{\tau} \mu$. But if $\lambda \xrightarrow{\tau} \mu$, then there exists integers $i, x, y$ with $i \geq 1$ and $0 \leq x < y$ such that
\begin{align*}
s(\lambda) &= s(\tau) \sqcup (x, x-e,\dotsc, x-(i-1)e), \\
s(\mu) &= s(\tau) \sqcup (y, y-e,\dotsc, y-(i-1)e),
\end{align*}
by Lemma \ref{L:ez2}(3).  Thus, $s(\lambda) < s(\mu)$ in the product order.
\end{proof}

\section{Setup}

In this section, we set up the notations which shall henceforth be used in this paper and state the main theorem.  We devote the next two sections to the proof of the main theorem.

We shall consider partitions with a fixed $e$-core, say $\kappa$.  We display all these partitions on an abacus with $l(\kappa) + Ne$ beads, where $N$ is `sufficiently large'.  For each $0 \leq i < e$, let $n_i$ be the number of beads on runner $i$ of such an abacus.  Define $\prec$ ($= \prec_{\kappa})$ as follows: $i \preceq j$ if and only if any of the following holds:
\begin{itemize}
\item $n_i < n_j$, or
\item $n_i = n_j$ and $i \leq j$.
\end{itemize}
We note that $\preceq$ is independent of $N$, and is a partial order on $\{0,1,\dotsc,e-1\}$.  Furthermore, as the abacus has $l(\kappa)+Ne$ beads, we have $n_0 \leq n_i$ for all $0 < i < e$, so that $0$ is the minimal element with respect to $\preceq$; we denote the maximal element by $M$.

For $i \prec M$ (which means $i \preceq M$ and $i \ne M$), define $i^+$ as the least (with respect to $\preceq$) $k$ such that $i \prec k$, and for $0 \prec j$, define $j^-$ as the largest (with respect to $\preceq$) $k$ such that $k \prec j$.

For $0 < i < e$, define
$$
d_i =
\begin{cases}
n_i - n_{i^-}, &\text{if } i^- < i;\\
n_i - n_{i^-} - 1, &\text{if } i^- > i.
\end{cases}
$$

Let $\pi$ be defined recursively as follows: $\pi(0) = 0$ and $\pi(i^+) = \pi(i) + 1$ for all $0 \prec i \prec M$.  Thus, $\pi$ is permutation on the set $\{0,1,\dotsc,e-1\}$, and $\pi(i) \leq \pi(j)$ if and only if $i \preceq j$.

When there is a need to indicate the $e$-core $\kappa$ for which $\preceq$, $M$, $i^{\pm}$, $d_i$ and $\pi$ are defined, we will include $\kappa$ as a subscript or superscript; we thus have $\preceq_{\kappa}$, $M_{\kappa}$, $i^{\pm_{\kappa}}$, $d_i^{\kappa}$ and $\pi_{\kappa}$.

The following Lemma is clear.

\begin{lem}
Let $\kappa$ be an $e$-core partition.  The following statements are equivalent:
\begin{enumerate}
\item $\kappa$ has an abacus display in which the number of beads in each runner is non-decreasing as we go from left to right.
\item $i \preceq j$ if and only if $i \leq j$.
\item $i^+ = i+1$ for all $0 \leq i \leq e-2$.
\item $i^- = i-1$ for all $0 < i < e$.
\item $\pi(i) = i$ for all $0 \leq i < e$.
\end{enumerate}

When any of these statements hold, we say $\kappa$ is of {\em Rouquier type}.
\end{lem}

We note that $\preceq_{\kappa}$, $M_{\kappa}$, $i^{\pm_{\kappa}}$, $d_i^{\kappa}$ and $\pi_{\kappa}$ and $\preceq_{\kappa'}$, $M_{\kappa'}$, $i^{\pm_{\kappa'}}$, $d_i^{\kappa'}$ and $\pi_{\kappa'}$ are related as follows:

\begin{lem} \label{L:kappa'}
Let $\kappa$ be an $e$-core partition.  Let $\Phi$ be the involution on $\{0,1,\dotsc,e-1\}$ defined by $\Phi(i) \equiv M_{\kappa} - i \pmod e$.  Then
\begin{enumerate}
\item $i \preceq_{\kappa} j$ if and only if $\Phi(j) \preceq_{\kappa'} \Phi(i)$;
\item $M_{\kappa} = M_{\kappa'}$;
\item $\Phi(i^{\pm_{\kappa'}}) = \Phi(i)^{\mp_{\kappa}}$, $\Phi(i^{\pm_{\kappa}}) = \Phi(i)^{\mp_{\kappa'}}$;
\item $d_i^{\kappa'} = d_{\Phi(i)^{+_\kappa}}^\kappa$;
\item $\pi_{\kappa'}(i) = e-1 - \pi_{\kappa}(\Phi(j))$.
\end{enumerate}
\end{lem}

\begin{proof}
This follows from the fact that an abacus display of $\kappa'$ may be obtained from that of $\kappa$ by rotating it through an angle of $\pi$ and reading the occupied positions as empty and the empty positions as occupied.
\end{proof}

Given two partitions $\lambda$ and $\mu$ with the same $e$-core, say $\kappa$, and $e$-quotients $(\lambda^0,\dotsc,\lambda^{e-1})$ and $(\mu^0,\dotsc, \mu^{e-1})$ respectively, define
\begin{align*}
\delta(\lambda,\mu) &= \sum_{j=1}^{e-1} \pi(j)(|\mu^j|-|\lambda^j|);\\
C_{\lambda\mu} &= \sum\prod_{j=0}^{e-1} c^{\lambda^j}_{\alpha^j\beta^j} c^{\mu^j}_{(\beta^{j^-})'\alpha^j},
\end{align*}
where the second sum runs over all partitions $\alpha^0, \dotsc, \alpha^{e-1}$ and $\beta^0,\dotsc,\beta^{M-1},\beta^{M+1},\dotsc,\beta^{e-1}$, and where $\beta^{M} = \emptyset = \beta^{0^-}$.
For convenience, we define $C_{\lambda\mu}$ and $\delta(\lambda,\mu)$ to be $0$ if $\lambda$ and $\mu$ do not have the same $e$-core. %Furthermore, write
%$$H(\mu) = \sum_{\lambda \in \mathcal{P}} C_{\lambda\mu} v^{\delta(\lambda,\mu)} \lambda.$$

We note the following:

\begin{lem} \label{beadgap}
Let $\mu $ and $\lambda $ be partitions with $e$-core $\kappa$, and $e$-quotients $(\mu^0,\mu^1,\dotsc,\mu^{e-1})$
and $(\lambda ^{0},\lambda ^{1},\dotsc ,\lambda ^{e-1})$ respectively, such
that for some partitions $\alpha^0,\dotsc ,\alpha^{e-1}$ and $\beta^0,\dotsc,\beta^{M-1},\beta^{M+1},\dotsc
,\beta^{e-1}$,
$$
\prod_{i=0}^{e-1} \left(c^{\lambda^i}_{\alpha_i\beta_i} c^{\mu^i}_{(\beta^{i^-})'\alpha^i} \right) \neq 0,
$$
where $\beta^{0^-}=\beta^{M}=\emptyset $.  Then
\begin{enumerate}
\item $\delta(\lambda,\mu) = \sum_{0\preceq i \prec M} |\beta^i|$;
\item whenever $x$ is the least vacant position on runner $j$, $y$ is the largest occupied position on runner $i$, and $i \prec j$, we have
$$
x-y >
\begin{cases}
(d_j+1-|\mu^{j^-}|-|\mu^j|-|\mu^{j^+}|)e, &\text{if } i^+ = j; \\
(\sum_{i \prec r \preceq j} d_r  + 1 - |\mu^i|-|\mu^{i^+}| - |\mu^j| - |\mu^{j^+}|)e, &\text{if } i^+ \prec j.
\end{cases}
$$
\end{enumerate}
\end{lem}

\begin{proof} \hfill
\begin{enumerate}
\item Since $\prod_{i=0}^{e-1} \left(c^{\lambda^i}_{\alpha_i\beta_i} c^{\mu^i}_{(\beta^{i^-})'\alpha^i} \right) \neq 0$, we have, for all $i$, $|\lambda^i| = |\alpha_i| + |\beta_i|$ and $|\mu^i| = |\beta^{i^-}| + |\alpha^i|$, so that $|\lambda^i| - |\mu^i| = |\beta^i| - |\beta^{i^-}|$.  Thus,
\begin{align*}
\delta(\lambda,\mu) &= \sum_{j=1}^{e-1} \pi(j) (|\mu^j| - |\lambda^j|) \\
&= \sum_{j=1}^{e-1} \pi(j) (|\beta^{j^-}| - |\beta^j|) \\
&= \sum_{j=1}^{e-1} \pi(j)|\beta^{j^-}| - \sum_{j=1}^{e-1} \pi(j)|\beta^j| \\
&= \sum_{0 \preceq j \prec M} (\pi(j)+1) |\beta^j| - \sum_{0 \prec j \preceq M} \pi(j)|\beta^j| \\
&= \sum_{0 \preceq j \prec M} |\beta^j|.
\end{align*}

\item Let $n_k$ be the number of  beads on runner $k$ in the abacus display of $\kappa$.  Then $x = (n_j - l(\lambda^j))e + j$, and $y = (n_i + \lambda^i_1 -1)e + i$.  Suppose $i \prec j$.  Then
\begin{align*}
x-y &= (n_j - n_i + 1 - l(\lambda^j) - \lambda^i_1)e + (j-i) \\
&= \left(\sum_{i \prec r \preceq j} (n_r - n_{r^-}) + 1 - l(\lambda^j) - \lambda^i_1\right)e + (j-i)
\end{align*}
If $i < j$, then clearly
\begin{align*}
x-y &\geq (\sum_{i \prec r \preceq j} d_r + 1 - l(\lambda^j) - \lambda^i_1)e + (j-i) \\
&> (\sum_{i \prec r \preceq j} d_r + 1 - l(\lambda^j) - \lambda^i_1)e.
\end{align*}
On the other hand, if $i > j$, then there exists $i \prec r \preceq j$ such that $r^- > r$, so that $n_r - n_{r^-} = d_r + 1$; thus
\begin{align*}
x-y &\geq (\sum_{i \prec r \preceq j} d_r + 2 - l(\lambda^j) - \lambda^i_1)e + (j-i) \\
&> (\sum_{i \prec r \preceq j} d_r + 1 - l(\lambda^j) - \lambda^i_1)e.
\end{align*}
Part (2) now follows since
$$l(\lambda^j) + \lambda^i_1 \leq |\lambda^j| + |\lambda^i| \leq
\begin{cases}
|\mu^{j^-}| + |\mu^j| + |\mu^{j^+}|, &\text{if } i^+=j;\\
|\mu^i|+|\mu^{i^+}| + |\mu^j| + |\mu^{j^+}|, &\text{if } i^+ \prec j.
\end{cases}
$$
\end{enumerate}
\end{proof}

\begin{Def}
Let $\mathcal{P}_{\kappa}^*$ be the collection of partitions $\mu$ having $e$-core $\kappa$ and $e$-quotient $(\mu^0,\dotsc, \mu^{e-1})$ such that
\begin{enumerate}
\item $|\mu^{i^-}| + |\mu^i| + |\mu^{i^+}| \leq d_i + 1$ for all $i=1,\dotsc, e-1$,
\item whenever $i$ and $j$ ($i \prec j$) satisfy
\begin{enumerate}
\item $|\mu^j| + |\mu^{j^+}| = d_j + 1$,
\item $|\mu^{i^-}| + |\mu^i| = d_i + 1$,
\end{enumerate}
there exists $k$ such that $i \prec k \prec j$ and $d_k > 0$.
\end{enumerate}
Here, $|\mu^{i^-}|$ (resp.\ $|\mu^{i^+}|$) is to be read as 0 when $i^-$ (resp.\ $i^+$) is undefined.
\end{Def}

\begin{eg*}
Let $e = 5$ and $\kappa = (3,3)$.  Then $(d_1, d_2, d_3, d_4) = (0,0,1,0)$.  Let $\mu = (8,3,2,1,1,1)$.  Then $\mu$ has $5$-core $\kappa$ and $5$-quotient $(\emptyset,(1),\emptyset,\emptyset,(1))$, so that $\mu \in \mathcal{P}_{\kappa}^*$ (but $\mu$ is not an element of $\mathcal{P}_{\kappa}$ defined in \cite{can}).
\end{eg*}

\begin{rem} \hfill
\begin{enumerate}
\item As a partition having $e$-core $\kappa$ is Rouquier if and only if its $e$-weight is not more than $\min\{d_i : 1 \leq i < e \} + 1$, we see that $\mathcal{P}_{\kappa}^*$ includes all Rouquier partitions with $e$-core $\kappa$.
\item In \cite{can}, a collection $\mathcal{P}_\kappa$ of partitions with $e$-core $\kappa$ is defined when $\kappa$ is of Rouquier type.  This is a subcollection of $\mathcal{P}_\kappa^*$.
\end{enumerate}
\end{rem}

The partitions in $\mathcal{P}_{\kappa}^*$ have the following nice property.

\begin{lem} \label{L:imptprop}
If $\mu \in \mathcal{P}_{\kappa}^*$, and $C_{\lambda\mu} \ne 0$, then any occupied position on runner $i$ of the abacus display of $\lambda$ is less than any vacant position on runner $j$ as long as $i \prec j$.
\end{lem}

\begin{proof}
This follows directly from the definition of $\mathcal{P}_{\kappa}^*$ and Lemma \ref{beadgap}.
\end{proof}

We now state the main theorem of this paper.

\begin{thm} \label{T:main}
Let $\mu \in \mathcal{P}_{\kappa}^*$ with $e$-quotient $(\mu^0,\dotsc, \mu^{e-1})$ .  Then
\begin{enumerate}
\item $d_{\lambda \mu}(v) = C_{\lambda\mu}v^{\delta(\lambda,\mu)}$ for all $\lambda \in \mathcal{P}$;
\item $d_{\lambda \mu}^l = C_{\lambda \mu}$ for all $\lambda \in \mathcal{P}$ and $l > \max_i (|\mu^i|)$.
\end{enumerate}
\end{thm}

We shall devote the next two sections to the proof of Theorem \ref{T:main}.

It is easy to describe the image under the Mullineux involuation of an $e$-regular partition for which Theorem \ref{T:main}(1) holds; furthermore, Theorem \ref{T:main}(1) also holds for the image:

\begin{prop} \label{P:Mull}
Suppose that an $e$-regular partition $\mu$ has $e$-core $\kappa$ and $e$-quotient $(\mu^0,\dotsc,\mu^{e-1})$, and that $d_{\lambda \mu}(v) = C_{\lambda\mu}v^{\delta(\lambda,\mu)}$ for all $\lambda \in \mathcal{P}$.  Then
\begin{enumerate}
\item $\mu^0 = \emptyset$;
\item $\mu^*$ has $e$-core $\kappa'$ and $e$-quotient $(\emptyset, \mu^{\Phi(1)^{+_\kappa}}, \mu^{\Phi(2)^{+_\kappa}}, \dotsc, \mu^{\Phi(e-1)^{+_\kappa}})$, where $\Phi$ is the involution on $\{ 0,1,\dotsc,e-1\}$ defined by $\Phi(i) \equiv M_{\kappa}-i \pmod e$;
\item $d_{\tau \mu^*}(v) = C_{\tau\mu^*}v^{\delta(\tau,\mu^*)}$ for all $\tau \in \mathcal{P}$.
\end{enumerate}
\end{prop}

\begin{proof}
By Theorem \ref{T:Mull2}, ${\mu^*}'$ is the unique partition such that $d_{{\mu^*}'\mu}(v)$ has degree equals the $e$-weight of $\mu$.  When $\prod_{j=0}^{e-1} c^{\lambda^j}_{\alpha^j\beta^j} c^{\mu^j}_{(\beta^{j^{-_{\kappa}}})'\alpha^j} \ne 0$ with $\beta^{0^-} = \beta^{M_{\kappa}} = \emptyset$, we have $\deg(d_{\lambda\mu}(v)) = \delta(\lambda,\mu) = \sum_{0 \preceq j \prec M_{\kappa}} |\beta^j|$ by Lemma \ref{beadgap}(1).  Thus $\deg(d_{\lambda\mu}(v)) = \sum_{i=0}^{e-1} |\mu^i|$ only if $(\beta^{j^{-_{\kappa}}})' = \mu^j$ and $\alpha^j = \emptyset$ for all $0 \preceq j \prec M_{\kappa}$.  It follows that $\mu^0 = \beta^{0^-} = \emptyset$, giving part (1).

Furthermore, the $e$-quotient of ${\mu^*}'$ equals
$$(\beta^0,\beta^1,\dotsc,\beta^{e-1}) = ((\mu^{0^{+_\kappa}})', (\mu^{1^{+_{\kappa}}})',\dotsc, (\mu^{(e-1)^{+_\kappa}})'),$$
where $\mu^{M^{+_\kappa}}$ is to be read as $\emptyset$.

Now, when a partition has $e$-core $\kappa$, and $e$-quotient $(\lambda^0,\dotsc,\lambda^{e-1})$ say, then its conjugate has $e$-core $\kappa'$ and $e$-quotient $((\lambda^{\Phi(0)})', (\lambda^{\Phi(1)})', \dotsc, (\lambda^{\Phi(e-1)})')$, where $\Phi$ is the involution on $\{ 0,1,\dotsc,e-1\}$ defined by $\Phi(i) \equiv M_{\kappa}-i \pmod e$.  Thus the $e$-quotient of $\mu^*$ is
$$(\mu^{\Phi(0)^{+_\kappa}}, \mu^{\Phi(1)^{+_\kappa}}, \mu^{\Phi(2)^{+_\kappa}}, \dotsc, \mu^{\Phi(e-1)^{+_\kappa}}) = (\emptyset, \mu^{\Phi(1)^{+_\kappa}}, \mu^{\Phi(2)^{+_\kappa}}, \dotsc, \mu^{\Phi(e-1)^{+_\kappa}}),$$ giving part (2).

Using Theorem \ref{T:Mull}, we need to show that $C_{\lambda'\mu^*} = C_{\lambda\mu}$ and $\delta(\lambda',\mu^*) = \sum_{j=1}^{e-1} |\mu_j| - \delta(\lambda,\mu)$ for part (3).  Since $\mu^0 = \emptyset$, we can simplify the expression of $C_{\lambda\mu}$, where $\lambda$ has $e$-core $\kappa$ and $e$-quotient $(\lambda^0,\dotsc,\lambda^{e-1})$, to
$$
C_{\lambda\mu} = \sum \left( \prod_{\substack{0 \leq i \leq e-1 \\[3pt] i \ne 0,\: M_\kappa}} c^{\lambda^i}_{\alpha^i\beta^i}
\prod_{\substack{0 \leq j \leq e-1 \\[3pt] j \ne 0,\: 0^{+_\kappa},\: M_\kappa }} c^{\mu^j}_{(\beta^{j^{-_\kappa}})'\alpha^j} \right)
c^{\mu^{0^{+_\kappa}}}_{(\lambda^0)'\alpha^{0^{+_\kappa}}} c^{\mu^{M_\kappa}}_{(\beta^{M_\kappa^{-_\kappa}})'\lambda^{M_\kappa}},
$$
where the sum runs over all partitions $\alpha^1,\dotsc,\alpha^{M_\kappa-1}, \alpha^{M_\kappa+1},\dotsc, \alpha^{e-2}$ and $\beta^1,\dotsc,\beta^{M_\kappa-1}, \beta^{M_\kappa+1},\dotsc, \beta^{e-2}$.  Similarly, we have
$$
C_{\lambda'\mu^*} = \sum \left( \prod_{\substack{0 \leq i \leq e-1 \\[3pt] i \ne 0,\: M_{\kappa'}}} c^{\nu^i}_{\gamma^i\delta^i} \prod_{\substack{0 \leq j \leq e-1 \\[3pt] j \ne 0,\: 0^{+_{\kappa'}},\: M_{\kappa'} }} c^{\rho^j}_{(\delta^{j^{-_{\kappa'}}})'\gamma^j} \right)
c^{\rho^{0^{+_{\kappa'}}}}_{(\nu^0)'\gamma^{0^{+_{\kappa'}}}} c^{\rho^{M_{\kappa'}}}_{(\delta^{M_{\kappa'}^{-_{\kappa'}}})'\nu^{M_{\kappa'}}},
$$
where the sum runs over all partitions $\gamma^1,\dotsc,\gamma^{M_{\kappa'}-1}, \gamma^{M_{\kappa'}+1},\dotsc, \gamma^{e-2}$ and $\delta^1,\dotsc,\delta^{M_{\kappa'}-1}, \delta^{M_{\kappa'}+1},\dotsc, \delta^{e-2}$, and $e$-quotients of $\lambda'$ and $\mu^*$ are denoted as $(\nu^0,\dotsc,\nu^{e-1})$ and $(\emptyset, \rho^1,\dotsc,\rho^{e-1})$ respectively.  We now simplify $C_{\lambda'\mu^*}$ using the fact that $\nu^i = (\lambda^{\Phi(i)})'$ and $\rho^i = \mu^{\Phi(i)^{+_\kappa}}$ for all $i$, and Lemma \ref{L:kappa'}; we have
\begin{align*}
C_{\lambda'\mu^*} &= \sum \left( \prod_{\substack{0 \leq i \leq e-1 \\[3pt] i \ne 0,\: M_{\kappa}}} c^{(\lambda^{\Phi(i)})'}_{\gamma^i\delta^i} \prod_{\substack{0 \leq j \leq e-1 \\[3pt] j \ne 0,\: 0^{+_{\kappa'}},\: M_{\kappa} }} c^{\mu^{\Phi(j^{-_{\kappa'}})}}_{(\delta^{j^{-_{\kappa'}}})'\gamma^j} \right)
c^{\mu^{M_{\kappa}}}_{\lambda^{M_{\kappa}} \gamma^{0^{+_{\kappa'}}}} c^{\mu^{0^+}}_{(\delta^{M_{\kappa}^{-_{\kappa'}}})'(\lambda^{0})'} \\
&= \sum \left( \prod_{\substack{0 \leq i \leq e-1 \\[3pt] i \ne 0,\: M_{\kappa}}} c^{(\lambda^{\Phi(i)})'}_{\gamma^i\delta^i} \prod_{\substack{0 \leq j \leq e-1 \\[3pt] j \ne 0,\: M_{\kappa}^{-_{\kappa'}},\: M_{\kappa} }} c^{\mu^{\Phi(j)}}_{(\delta^{j})'\gamma^{j^{+_{\kappa'}}}} \right)
c^{\mu^{M_{\kappa}}}_{\lambda^{M_{\kappa}} \gamma^{0^{+_{\kappa'}}}} c^{\mu^{0^{+_\kappa}}}_{(\delta^{M_{\kappa}^{-_{\kappa'}}})'(\lambda^{0})'} \\
&= \sum \left( \prod_{\substack{0 \leq i \leq e-1 \\[3pt] i \ne 0,\: M_{\kappa}}} c^{(\lambda^i)'}_{\gamma^{\Phi(i)}\delta^{\Phi(i)}} \prod_{\substack{0 \leq j \leq e-1 \\[3pt] j \ne 0,\: 0^{+_\kappa},\: M_{\kappa}}} c^{\mu^j}_{(\delta^{\Phi(j)})'\gamma^{\Phi(j^{-_\kappa})}} \right)
c^{\mu^{M_{\kappa}}}_{\lambda^{M_{\kappa}} \gamma^{\Phi(M_{\kappa}^{-_\kappa})}} c^{\mu^{0^{+_\kappa}}}_{(\delta^{\Phi(0^{+_ \kappa})})'(\lambda^{0})'} \\
&=C_{\lambda\mu}.
\end{align*}
When $C_{\lambda\mu}\ (= C_{\lambda'\mu^*}) \ne 0$, then $d_{\lambda\mu}(v) \ne 0$, so that $\lambda$ and $\mu$ must have the same $e$-weight, i.e.\ $\sum_{i=0}^{e-1} |\lambda^i| = \sum_{i=1}^{e-1} |\mu^i|$ (since $|\mu^0| = 0$).  Using Lemma \ref{L:kappa'}, we have
{\allowdisplaybreaks
\begin{align*}
\delta(\lambda',\mu^*)
&= \sum_{j=1}^{e-1} \pi_{\kappa'}(j)(|\mu^{\Phi(j)^{+_\kappa}}| - |(\lambda^{\Phi(j)})'|)\\
&= \sum_{j=1}^{e-1} (e-1-\pi_\kappa(\Phi(j)))(|\mu^{\Phi(j)^{+_\kappa}}| - |\lambda^{\Phi(j)}|)\\
&= (e-1) \sum_{j=1}^{e-1} |\mu^{\Phi(j)^{+_\kappa}}| - (e-1) \sum_{j=1}^{e-1} |\lambda^{\Phi(j)}| - \sum_{j=1}^{e-1} \pi (\Phi(j))(|\mu^{\Phi(j)^{+_\kappa}}| - |\lambda^{\Phi(j)}|) \\
&= (e-1) |\lambda^{M_\kappa}| - \sum_{\substack{0 \leq j \leq e-1 \\[3pt] j \ne M_\kappa}} \pi(j)(|\mu^{j^{+_\kappa}}| - |\lambda^j|) \\
&= (e-1) |\lambda^{M_\kappa}| - \sum_{\substack{0 \leq j \leq e-1 \\[3pt] j \ne M_\kappa}} \pi(j)|\mu^{j^{+_\kappa}}| + \sum_{\substack{0 \leq j \leq e-1 \\[3pt] j \ne M_\kappa}} \pi(j)|\lambda^j| \\
&= (e-1) |\lambda^{M_\kappa}| - \sum_{j=1}^{e-1} (\pi(j)-1)|\mu^{j}| + \sum_{\substack{0 \leq j \leq e-1 \\[3pt] j \ne M_\kappa}} \pi(j)|\lambda^j| \\
&= \sum_{j=1}^{e-1} |\mu^j| - \sum_{j=1}^{e-1} \pi(j)(|\mu^j| - |\lambda^{j}|) \\
&= \sum_{j=1}^{e-1} |\mu^j| - \delta(\lambda,\mu).
\end{align*}}
This completes the proof.
\end{proof}

We end this section by obtaining some corollaries on the $e$-regular partitions contained in $\mathcal{P}_\kappa^*$ while assuming Theorem \ref{T:main}.

Firstly, the $e$-regular partitions in $\mathcal{P}_{\kappa}^*$ can be easily described:

\begin{cor}
Let $\mu \in \mathcal{P}_{\kappa}^*$, with $e$-quotient $(\mu^0,\dotsc,\mu^{e-1})$.  Then $\mu$ is $e$-regular if and only if $\mu^0 = \emptyset$.
\end{cor}

\begin{proof}
Suppose $\mu$ is $e$-regular.  Then $\mu^0 = \emptyset$ by Theorem \ref{T:main}(1) and Proposition \ref{P:Mull}(1).  Conversely, suppose $\mu^0 = \emptyset$.  Then the least $N$ positions on runner $0$ in the abacus display of $\mu$ (with $l(\kappa) + Ne$ beads) are occupied, while the rest are vacant.  By Lemma \ref{L:imptprop} (with $i = 0$), the least $(N-1)$ positions on runner $j$ are all occupied for all $j >0$.  It follows that $\mu$ is $e$-regular.
\end{proof}

Also, $\mathcal{P}_{\kappa}^*$ is closed under the Mullineux involution:

\begin{cor}
Suppose $\mu \in \mathcal{P}_{\kappa}^*$ is $e$-regular.  Then $\mu^* \in \mathcal{P}_{\kappa'}^*$.
\end{cor}

\begin{proof}
This follows directly from Proposition \ref{P:Mull}(2), Lemma \ref{L:kappa'} and the definitions of
$\mathcal{P}_{\kappa}^*$ and $\mathcal{P}_{\kappa'}^*$.
\end{proof}

\section{Canonical basis}

In this section, we provide the proof of part (1) of Theorem \ref{T:main}.

We begin with a reformulation of the Theorem.

\begin{Def}
For any partition $\mu$, let $H(\mu) = \sum_{\lambda \in \mathcal{P}} C_{\lambda\mu} v^{\delta(\lambda,\mu)} \lambda$.
\end{Def}

Theorem \ref{T:main}(1) can thus be reformulated as follows:

If $\mu \in \mathcal{P}_{\kappa}^*$, then $G(\mu) = H(\mu)$.

\begin{lem} \label{L:barinvariant}
$G(\mu) = H(\mu)$ if and only if $\overline{H(\mu)} = H(\mu)$.
\end{lem}

\begin{proof}
Note that $\delta(\lambda,\mu)$ is always non-negative.  Furthermore, $C_{\lambda\mu} \ne 0$ and $\delta(\lambda,\mu) = 0$ if and only if $\lambda = \mu$, and $C_{\mu\mu} = 1$.  Thus, $H(\mu) - \mu \in \bigoplus_{\lambda} v\mathbb{Z}[v] \lambda$ and the lemma follows.
\end{proof}

We now review the results of \cite{can}.
For an $e$-core partition $\kappa$ of Rouquier type, and integers $a$ and $k$ with $1 \leq a \leq e-1$ and $k\geq 1$, let
$$
F_{a,k} = \mathfrak{f}_a^{(k)} \mathfrak{f}_{a+1}^{(k)} \dotsm \mathfrak{f}_{e-1}^{(k)} \mathfrak{f}_{a-1}^{(k)} \mathfrak{f}_{a-2}^{(k)} \dotsm \mathfrak{f}_0^{(k)},$$
where, if $j$ is the residue class of $i + l(\kappa)$ modulo $e$, then $\mathfrak{f}_j = f_i \in U_v(\widehat{\mathfrak{sl}}_e)$.  Note that $\overline{F_{a,k}(x)} = F_{a,k}(\overline{x})$ for all $x \in \mathcal{F}$.

\begin{lem}[{\cite[Lemma 3.1]{can}}] \label{L:messy}
Let $\kappa$ be an $e$-core partition of Rouquier type.  Let $\lambda$ be a partition with $\kappa$ and $e$-quotient $(\lambda^0,\dotsc, \lambda^{e-1})$, and consider its abacus display.  Let $x_a$ and $y_a$ (resp.\ $x_{a-1}$ and $y_{a-1}$) be respectively the least vacant position and largest occupied position on runner $a$ (resp.\ $a-1$).  Suppose that
\begin{enumerate}
\item[(i)] $x_{a-1} - u > e$ for all occupied positions $u$ lying on a runner to the left of runner $a-1$;
\item[(ii)] $t - y_a > e$ for all vacant positions $t$ lying on a runner to the right of runner $a$.
\end{enumerate}
If $1 \leq k \leq (x_a - y_{a-1} -1)/e$, then
$$F_{a,k} (\lambda) = \sum_{j=0}^k v^j \sum_{\alpha,\beta} c^{\alpha}_{\lambda^{a-1}(j)} c^{\beta}_{\lambda^a(1^{k-j})}\, \lambda(\alpha,\beta),$$
where $\lambda(\alpha,\beta)$ denotes the partition having $e$-core $\kappa$ and $e$-quotient
$$(\lambda^0, \dotsc, \lambda^{a-2}, \alpha, \beta, \lambda^{a+1},\dotsc, \lambda^{e-1}).$$
\end{lem}

\begin{note}
If $(\prod_{j=0}^{e-1} c^{\lambda^j}_{\alpha^j\beta^j} c^{\mu^j}_{(\beta^{j-1})'\alpha^j}) \ne 0$ with $\beta^{-1} = \beta^{e-1} = \emptyset$, then $\delta(\lambda,\mu) = \sum_{j=0}^{e-2} |\beta^j|$.
\end{note}

The proof of Proposition 4.1 of \cite{can} can be easily adapted to prove the following proposition.

\begin{prop} \label{P:induct}
Let $\kappa$ be an $e$-core partition of Rouquier type.
Let $\mu$ be a partition with $e$-core $\kappa$ and $e$-quotient $(\mu^0,\dotsc, \mu^{e-1})$, and suppose that for all $\lambda$ such that $C_{\lambda\mu} \ne 0$, we have
$$F_{a,k} (\lambda) = \sum_{j=0}^k v^j \sum_{\alpha,\beta} c^{\alpha}_{\lambda^{a-1}(j)} c^{\beta}_{\lambda^a(1^{k-j})}\, \lambda(\alpha,\beta).$$
Then
$$F_{a,k} (H(\mu)) = \sum_{\eta} c^{\eta}_{\mu^a(1^k)} H(\mu_\eta),$$ where $\mu_\eta$ denotes the partition having $e$-core $\kappa$ and $e$-quotient
$$(\mu^0,\dotsc,\mu^{a-1}, \eta, \mu^{a+1},\dotsc \mu^{e-1}).$$
\end{prop}

The following is a direct consequence arising from the definition of $\mathcal{P}_{\kappa}^*$ and Lemma \ref{beadgap}.

\begin{lem} \label{L:ez}
Let $\kappa$ be an $e$-core partition of Rouquier type.
Let $\mu \in \mathcal{P}_{\kappa}^*$ with $e$-quotient $(\mu^0,\dotsc,\mu^{e-1})$.  Suppose $\mu^a \ne \emptyset$ for some $a \geq 1$.  Let $\tau^a$ be any partition such that $|\tau^a| < |\mu^a|$, let $\tau$ be the partition having $e$-core $\kappa$ and $e$-quotient $(\mu^0,\dotsc, \mu^{a-1},\tau^a,\mu^{a+1},\dotsc,\mu^{e-1})$.  Then any $\lambda$ with $C_{\lambda\tau} \ne 0$, $a$ and $k = |\mu^a| - |\tau^a|$ satisfy the hypothesis of Lemma \ref{L:messy}.
\end{lem}

\begin{prop} \label{P:Rouq}
Let $\kappa$ be an $e$-core partition of Rouquier type.
Let $\mu \in \mathcal{P}_{\kappa}^*$ .  Then $G(\mu) = H(\mu)$.
\end{prop}

\begin{proof}
We prove by induction.  Clearly, $G(\kappa) = \kappa = H(\kappa)$.  Let $|\mu| > |\kappa|$, and suppose $G(\tau) = H(\tau)$ holds for all $\tau \in \mathcal{P}_{\kappa}^*$ satisfying either $|\tau| < |\mu|$ or $\tau < \mu$.  Let $(\mu^0,\dotsc, \mu^{e-1})$ be the $e$-quotient of $\mu$.

Suppose first that there exists $a\geq 1$ such that $\mu^a \ne \emptyset$.  Let $\tau^a$ be the partition obtained from $\mu^a$ by removing its first column, and let $\tau$ denote the partition having $e$-core $\kappa$ and $e$-quotient $(\mu^0,\dotsc, \mu^{a-1}, \tau^a,\mu^{a+1},\dotsc, \mu^{e-1})$.  Let $k = |\mu^a| - |\tau^a|$.  By Lemmas \ref{L:messy} and \ref{L:ez} and Proposition \ref{P:induct}, we have $F_{a,k} (H(\tau)) = \sum_{\eta} c^{\eta}_{\tau^a(1^k)} H(\tau_{\eta})$.  Now, $\tau, \tau_{\eta} \in \mathcal{P}^*_\kappa$ for all $\eta$ such that $c^{\eta}_{\tau^a(1^k)} \ne 0$.  Furthermore, $H(\mu)$ occurs exactly once as a summand of $F_{a,k} (H(\tau))$ and all the other summands $H(\tau_{\eta})$ satisfy $\tau_\eta < \mu$.  Thus, $F_{a,k} (H(\tau))$, and all its summands $H(\tau_\eta)$ with $\tau_\eta \ne \mu$ are bar-invariant by induction hypothesis, and hence so is $H(\mu)$.  Thus $G(\mu) = H(\mu)$ by Lemma \ref{L:barinvariant}.

It remains to consider the case where $\mu^i = \emptyset$ for all $i > 0$.  In this case, $H(\mu) = \mu$.  Let $\lambda$ be a partition having $e$-core $\kappa$ and $e$-quotient $(\lambda^0,\dotsc, \lambda^{e-1})$.  It is not difficult to see that if $\lambda^i \ne \emptyset$ for some $i>0$, then $\lambda > \mu$, so that $d_{\lambda\mu}(v) = 0$.  If $\lambda^i = \emptyset$ for all $i>0$, then we shall show later in Corollary \ref{C:=0} that $d_{\lambda\mu}^0 = 0$ if $\lambda \ne \mu$.  Since $d_{\lambda\mu}(1) = d_{\lambda\mu}^0$ and $d_{\lambda\mu}(v) \in \mathbb{N}_0[v]$, we thus have $d_{\lambda\mu}(v) = 0$ if $\lambda \ne \mu$.  Hence $G(\mu) = \mu = H(\mu)$.
\end{proof}

\begin{prop} \label{P:induct2}
Let $\kappa$ be an $e$-core partition, and for $0 \leq j < e$, let $n_j$ be the number of beads on runner $j$ of its abacus display.  Suppose that $n_{i-1} > n_i$ for some $i \geq 2$.  Given a partition $\lambda$ with $e$-core $\kappa$, let $\Psi(\lambda)$ be the partition obtained from $\lambda$ by interchanging the $(i-1)$-th and $i$-th runners of its abacus display.  Let $\mu \in \mathcal{P}_{\kappa}^*$.  Then
\begin{enumerate}
\item $\Psi(\mu) \in \mathcal{P}_{\Psi(\kappa)}^*$;
\item $f_r^{(k)}H(\mu) = H(\Psi(\mu))$;
\item $e_r^{(k)}H(\Psi(\mu)) = H(\mu)$,
\end{enumerate}
where $r$ is the residue class of $l(\kappa)-i$ modulo $e$, and $k = n_{i-1}-n_i$.
\end{prop}

\begin{proof}
Part (1) follows from the definition of $\mathcal{P}_{\kappa}^*$ (and $\mathcal{P}_{\Psi(\kappa)}^*$).  For parts (2) and (3), it suffices to show the following:
\begin{itemize}
\item for all partitions $\lambda$ with $e$-core $\kappa$, $C_{\lambda\mu} = C_{\Psi(\lambda)\Psi(\mu)}$ and $\delta(\lambda,\mu) = \delta(\Psi(\lambda),\Psi(\mu))$;
\item whenever $C_{\tau\mu} \ne 0$, the largest occupied position on runner $i$ of the abacus display of $\tau$ is less than the least vacant position on runner $i-1$.
\end{itemize}
The first assertion follows from definition, while the second follows from Lemma \ref{L:imptprop}.
\end{proof}

We are now ready to prove Theorem \ref{T:main}(1).

\begin{proof}[Proof of Theorem \ref{T:main}(1)]
We prove by induction on $N_{\kappa} = |\{(i,j) \mid i< j,\ j \prec_{\kappa} i \}|$.  If $N_{\kappa} = 0$, then $\kappa$ is of Rouquier type, so that $G(\mu) = H(\mu)$ by Proposition \ref{P:Rouq}.  Assume thus $N_{\kappa} >0$.  Then there exists $2 \leq i < e$ such that $n_{i-1} > n_i$ (here, and hereafter in this proof, we keep the notations of Proposition \ref{P:induct2}).  By Proposition \ref{P:induct2}(1), $\Psi(\mu) \in \mathcal{P}_{\Psi(\kappa)}^*$.  Since $N_{\Psi(\kappa)} < N_\kappa$, we have $G(\Psi(\mu)) = H(\Psi(\mu))$ by induction hypothesis.  By Proposition \ref{P:induct2}(3), we have $H(\mu) = e_r^{(k)}H(\Psi(\mu)) = e_r^{(k)}G(\Psi(\mu))$, so that $H(\mu)$ is bar-invariant.  Hence $G(\mu) = H(\mu)$ by Lemma \ref{L:barinvariant}.
\end{proof}

\section{$q$-Schur algebras}

In this section, we provide the proof of part (1) of Theorem \ref{T:main}, which may be reformulated as follows: If $\mu \in \mathcal{P}_{\kappa}^*$ with $e$-quotient $(\mu^0,\mu^1,\dotsc,\mu^{e-1})$, then 
$
d_{\lambda\mu}^l = d_{\lambda\mu}^0
$
for all $l > \max_i(|\mu^i|)$.

\begin{Def}
Let $\approx$ be the equivalence relation defined on $\mathcal{P}$ as follows:  $\lambda \approx \mu$ if and only if $\lambda$ and $\mu$ are partitions having the same $e$-core, and $|\lambda^i| = |\mu^i|$ for all $i$, where $(\lambda^0,\lambda^1,\dotsc, \lambda^{e-1})$ and $(\mu^0,\mu^1,\dotsc, \mu^{e-1})$ are the $e$-quotients of $\lambda$ and $\mu$ respectively.
\end{Def}

\begin{Def}
Let $\kappa$ be an $e$-core partition of Rouquier type.  Let
$\bar{\mathcal{P}}_{\kappa}$ be a collection of partitions with $e$-core $\kappa$, having the following two properties:
\begin{enumerate}
\item On the abacus display of any partition in $\bar{\mathcal{P}}_{\kappa}$, any pair $(x,y)$, where $x$ is an occupied position on runner $i$ while $y$ is a vacant position on runner $j$, with $i < j$, satisfies $x < y$;
\item If $\lambda, \mu \in \bar{\mathcal{P}}_{\kappa}$, $\lambda \approx \mu$, with respective $e$-quotients $(\lambda^0,\dotsc,\lambda^{e-1})$ and $(\mu^0,\dotsc,\mu^{e-1})$, then every partition having $e$-core $\kappa$, and $e$-quotient of the form $(\mu^0,\dotsc,\mu^{i-1},\lambda^{i},\dotsc,\lambda^{e-1})$ ($1\leq i \leq e-1$) lies in $\bar{\mathcal{P}}_{\kappa}$.
\end{enumerate}
\end{Def}

\begin{lem} \label{L:*bar}
Let $\kappa$ be an $e$-core partition of Rouquier type.  Then $\mathcal{P}_{\kappa}^*$ satisfies the two conditions of $\bar{\mathcal{P}}_{\kappa}$.
\end{lem}

\begin{proof}
If $\mu \in \mathcal{P}_{\kappa}^*$, then $\mu$ satisfies condition (1) of $\bar{\mathcal{P}}_{\kappa}$ by Lemma \ref{L:imptprop}.  That condition (2) of $\bar{\mathcal{P}}_{\kappa}$ also holds follows directly from the definition of $\mathcal{P}_\kappa^*$ .
\end{proof}

The following is an immediate consequence.

\begin{lem} \label{L:seq}
Let $\kappa$ be an $e$-core partition of Rouquier type.  Let $\lambda \in \bar{\mathcal{P}}_{\kappa}$, and let $s(\lambda) = (c_1,c_2,\dotsc, c_w)$ be its induced $e$-sequence.  Then $\overline{c_i} \geq \overline{c_{i+1}}$ for all $1 \leq i \leq w-1$, where $\overline{x}$ denotes the residue class of $x$ modulo $e$.
\end{lem}

\begin{proof}
This follows from Lemma \ref{L:ez2}(2) and condition (1) of $\bar{\mathcal{P}}_\kappa$.
\end{proof}

\begin{prop} \label{P:nleq_p}
Let $\kappa$ be an $e$-core partition of Rouquier type.  Suppose $\lambda, \mu \in \bar{\mathcal{P}}_{\kappa}$ with $\lambda \approx \mu$.  If $\lambda \to \tau$ with $\lambda \not\approx \tau$, then $\tau \nleq_p \mu$.
\end{prop}

\begin{proof}
Suppose $\lambda \xrightarrow{\sigma} \tau$.  By Lemma \ref{L:ez2}(3), there exist integers $N, x, y$ with $N \geq 1$ and $0 \leq x < y$ such that
\begin{align*}
s(\lambda) &= s(\sigma) \sqcup (x,x-e,\dotsc,x-(N-1)e),\\
s(\tau) &= s(\sigma) \sqcup (y,y-e,\dotsc,y-(N-1)e).
\end{align*}
Suppose $\overline{x} = i$ and $\overline{y} = j$, where $\overline{a}$ denotes the residue class of $a$ modulo $e$.  Condition (1) of $\bar{\mathcal{P}}_\kappa$ forces $j > i$.  Let $s(\lambda) = (l_1,l_2,\dotsc,l_w)$, and let $r$ be the least index such that $\overline{l_r} = i$.  Then $\overline{l_s} > i$ for all $s < r$ by Lemma \ref{L:seq}, so that $l_1,l_2,\dotsc,l_{r-1}$ are terms of $s(\sigma)$.  Let $s(\mu) = (m_1,m_2,\dotsc,m_w)$.  Then since $\lambda \approx \mu$, we have $\overline{m_t} = \overline{l_t}$ for all $t$ by Lemma \ref{L:seq}; in particular, $\overline{m_s} > i$ for all $s < r$.  Let $\rho$ denote the partition having $e$-core $\kappa$ and $e$-quotient $(\mu^0,\dotsc,\mu^i, \lambda^{i+1},\dotsc,\lambda^{e-1})$.  Then $\rho \in \bar{\mathcal{P}}_\kappa$, and $s(\rho) = (l_1,\dotsc,l_{r-1}, m_r,\dotsc, m_w)$ by Lemma \ref{L:seq}.  In the abacus display of $\rho$, $m_r$ is the largest occupied position on runner $i$, so it is less than all the vacant positions on runners to the right of runner $i$ by condition (1) of $\bar{\mathcal{P}}$; in particular, $m_r < y$ since $y$ is a vacant position on runner $j$ of $\lambda$ (and $j > i$), and hence of $\rho$.  Thus, $m_r$ is less than $r$ terms of $s(\tau)$, namely, $l_1,l_2,\dotsc, l_{r-1}, y$.  This implies that $\mu \ngeq_p \tau$.
\end{proof}

\begin{cor} \label{C:=0}
Let $\kappa$ be an $e$-core partition of Rouquier type.  Let $\lambda, \mu \in \bar{\mathcal{P}}_\kappa$ with $\lambda \approx \mu$.  If $l= 0$ or $l$ is greater than the size of each constituent of the $e$-quotient of $\lambda$, then $d_{\lambda\mu}^l = \delta_{\lambda\mu}$.
\end{cor}

\begin{proof}
By Theorem \ref{T:Jantzen}, it suffices to prove that $J_{\lambda\mu} = 0$ when $\lambda \ne \mu$.
For each $\tau$ such that $\lambda \to \tau$, let $j_{\tau\mu}=\sum_{\sigma} (-1)^{l_{\lambda\sigma} + l_{\tau\sigma} + 1}$, where $\sigma$ runs over all partitions such that $\lambda \xrightarrow{\sigma} \tau$.
By Proposition \ref{P:nleq_p} and Lemmas \ref{L:Jimplyp} and \ref{L:J}, $d_{\tau\mu}^l = 0$ for all $\tau$ such that $\lambda \to \tau$ and $\lambda \not\approx \tau$.  Furthermore, the condition on $l$ implies that we always have $\nu_l(h_{\lambda\sigma}) = 0$.  Thus $J_{\lambda\mu} = \sum_\tau j_{\tau\mu} d_{\tau\mu}^l$, where $\tau$ runs over all partitions such that $\lambda \to \tau$ and $\lambda \approx \tau$.  Fix such a $\tau$; suppose $\tau$ can be obtained from $\lambda$ by first moving a bead at position $x$ on runner $j$ to position $x-Ne$ and then moving a bead at position $y-Ne$ on runner $j$ to $y$, with $x < y$.  Let $\rho_1$ (resp.\ $\rho_2$) be the partition obtained from $\lambda$ by moving the bead at position $x$ (resp.\ $y-Ne$) to $x-Ne$.  Then $\lambda \xrightarrow{\sigma} \tau$ if and only if $\sigma = \rho_1$ or $\rho_2$.  Furthermore,
$l_{\lambda\rho_1} + l_{\tau\rho_1}$ and $l_{\lambda\rho_2} + l_{\tau\rho_2}$ are of different parity.  Thus $j_{\tau\mu} = 0$.  Hence $J_{\lambda\mu} = 0$.
\end{proof}

We need the following results, which are analogous to each other, and which we believe are well-known, but we are not able to find an appropriate reference in the existing literature.

\begin{prop} \label{P:gen}
Suppose $P_{\mathbb{C}}^\mu \up{} = \bigoplus_{i=1}^n a_i P_{\mathbb{C}}^{\sigma_i}$ (where $(- \uparrow)$ is a composition of (divided powers) of $r$-induction functors), $d_{\lambda\mu}^0 = d_{\lambda\mu}^l$ for all $\lambda$, and $d_{\sigma_i\sigma_j}^l = \delta_{ij}$.  Then $P_{\mathbb{F}}^\mu \up{} = \bigoplus_{i=1}^n a_i P_{\mathbb{F}}^{\sigma_i}$, and $d_{\tau\sigma_i}^0 = d_{\tau\sigma_i}^l$ for all $\tau$ and all $i$.
\end{prop}

\begin{prop} \label{P:gen2}
Suppose $P_{\mathbb{C}}^\mu \down{} = \bigoplus_{i=1}^n a_i P_{\mathbb{C}}^{\sigma_i}$ (where $(- \downarrow)$ is a composition of (divided powers) of $r$-restriction functors), $d_{\lambda\mu}^0 = d_{\lambda\mu}^l$ for all $\lambda$, and $d_{\sigma_i\sigma_j}^l = \delta_{ij}$.  Then $P_{\mathbb{F}}^\mu \down{} = \bigoplus_{i=1}^n a_i P_{\mathbb{F}}^{\sigma_i}$, and $d_{\tau\sigma_i}^0 = d_{\tau\sigma_i}^l$ for all $\tau$ and all $i$.
\end{prop}

We provide a proof of Proposition \ref{P:gen} only; that of Proposition \ref{P:gen2} is entirely analogous.

\begin{proof}[Proof of Proposition \ref{P:gen}]
Note that, since $\Delta_\mathbb{F}^\mu$ and $\Delta_\mathbb{C}^\mu$ have the same branching rule and $d_{\lambda\mu}^0 = d_{\lambda\mu}^l$ for all $\lambda$, we have $[P_{\mathbb{F}}^\mu \up{} : \Delta_{\mathbb{F}}^{\tau}]= [P_{\mathbb{C}}^\mu \up{} : \Delta_{\mathbb{C}}^{\tau}] = \sum_{i=1}^n a_i d_{\tau\sigma_i}^0 \leq \sum_{i=1}^n a_i d_{\tau\sigma_i}^l$.  We may assume that $\sigma_i >_J \sigma_j$ implies that $i < j$.  We prove by induction on $r$ that $\bigoplus_{i=1}^r a_i P_{\mathbb{F}}^{\sigma_i}$ is a direct summand of $P_{\mathbb{F}}^{\mu} \up{}$.  This is clear for $r = 1$.  Suppose $\bigoplus_{i=1}^{r-1} a_i P_{\mathbb{F}}^{\sigma_i}$ is a direct summand of $P_{\mathbb{F}}^{\mu} \up{}$, so that
\begin{equation*}
\sum_{i=1}^{r-1} a_i d_{\tau\sigma_i}^l = [\bigoplus_{i=1}^{r-1} a_i P_{\mathbb{F}}^{\sigma_i} : \Delta_\mathbb{F}^\tau] \leq [P_{\mathbb{F}}^\mu \up{} : \Delta_{\mathbb{F}}^{\tau}] \leq \sum_{i=1}^n a_i d_{\tau\sigma_i}^l. \tag{$*$}
\end{equation*}
for all $\tau$.  If $\tau >_J \sigma_r$, then $\sigma_j \ngeq_J \tau$ for all $j \geq r$, so that $d_{\tau\sigma_j}^l = 0$, and hence we have equality in $(*)$.  This shows that if $\tau >_J \sigma_r$, then $P_\mathbb{F}^{\tau}$ is not a summand of $P_{\mathbb{F}}^{\mu} \up{}$.  Now, $[\bigoplus_{i=1}^{r-1} a_i P_{\mathbb{F}}^{\sigma_i} : \Delta_\mathbb{F}^{\sigma_r}] = \sum_{i=1}^{r-1} a_i d_{\sigma_r\sigma_i}^l = 0$, while $[P_{\mathbb{F}}^\mu \up{} : \Delta_\mathbb{F}^{\sigma_r}] = a_r$.  Thus, $a_r P_{\mathbb{F}}^{\sigma_r}$ is a direct summand of $P_{\mathbb{F}}^\mu \up{}$.  Hence, by induction, $\bigoplus_{i=1}^n a_i P_{\mathbb{F}}^{\sigma_i}$ is a direct summand of $P_{\mathbb{F}}^{\mu} \up{}$, and so $(*)$ holds when $r-1$ is replaced by $n$, and in fact with equality throughout.  Thus $P_{\mathbb{F}}^\mu \up{} = \bigoplus_{i=1}^n a_i P_{\mathbb{F}}^{\sigma_i}$, and $d_{\tau\sigma_i}^0 = d_{\tau\sigma_i}^l$ for all $\tau$ and all $i$.
\end{proof}

We are now ready to prove Theorem \ref{T:main}(2).

\begin{proof}[Proof of Theorem \ref{T:main}(2)]
We only need to show the Theorem holds when $l > \max_i(|\mu^i|)$ and $\lambda$ has $e$-core $\kappa$ and $e$-quotient $(\lambda^0,\dotsc,\lambda^{e-1})$ with $\sum_{i=0}^{e-1} |\lambda^i| = \sum_{i=0}^{e-1} |\mu^i|$.

Let $N_{\kappa} = |\{ (i,j) \mid i< j,\ j \prec_\kappa i \}|$.  

Suppose first that $N_{\kappa} = 0$.  Then $\kappa$ is of Rouquier type.
\begin{description}
\item[Case 1. $\mu^i= \emptyset$ for all $i \geq 1$]  If $\lambda^i \ne \emptyset$ for some $i \geq 1$, then $\lambda > \mu$ so that $d_{\lambda\mu}^l = 0 = C_{\lambda\mu}$.  If $\lambda^i = \emptyset$ for all $i \geq 1$, then $\lambda \approx \mu$, so that $d_{\lambda\mu}^l = \delta_{\lambda\mu} = C_{\lambda\mu}$ by Corollary \ref{C:=0}.
\item[Case 2. $\mu^a \ne \emptyset$ for some $a \geq 1$] From the proof of Proposition \ref{P:Rouq}, we can find a $\tau \in \mathcal{P}_\kappa^*$ and an integer $k$ such that $G(\mu)$ occurs exactly once as a summand of $F_{a,k} G(\tau)$.  Furthermore, if $G(\alpha)$ and $G(\beta)$ are summands of $F_{a,k} G(\tau)$, then $\alpha \approx \beta$, so that $d_{\alpha,\beta}^l = \delta_{\alpha\beta}$ by Corollary \ref{C:=0}.  We may assume that $d_{\sigma\tau}^l = d_{\sigma\tau}^0$ for all $\sigma \in \mathcal{P}$ by induction hypothesis.  Thus, $d_{\lambda\mu}^l = d_{\lambda\mu}^0$ for all $\lambda \in \mathcal{P}$ by Proposition \ref{P:gen}.
\end{description}

Now, suppose $N_{\kappa} > 0$.  Then there exists $2 \leq i < e$ such that $n_{i-1} > n_i$ (here, and hereafter in this proof, we keep the notations of Proposition \ref{P:induct2}).  By Proposition \ref{P:induct2}, we have $e_r^{(k)}G(\Psi(\mu)) = G(\mu)$.  Since $N_{\Psi(\kappa)} < N_{\kappa}$, we may assume that $d_{\sigma\Psi(\mu)}^0 = d_{\sigma\Psi(\mu)}^l$ for all $\sigma \in \mathcal{P}$ by induction hypothesis.  Thus $d_{\lambda\mu}^0 = d_{\lambda\mu}^l$ for all $\lambda \in \mathcal{P}$ by Proposition \ref{P:gen2}.
\end{proof}

 \end{document}